\input amstex
\documentstyle{amsppt}
\magnification=\magstep1
\define \my{\bf}
\define\po{\parindent 0pt}
\define \mytitle{\po \m \my}
\define \p{{\po \it \m Proof. }}
\define \m{\medskip}
\nologo
\pageheight{16.5cm}
\TagsOnLeft
\define \ph#1{\phantom {#1}}
\define \T#1{\widetilde {#1}}
\define \Ker{\operatorname{Ker}} 
\define \IM{\operatorname {Im}}
\define \cat{\operatorname{cat}}
\define \swgt{\operatorname{swgt}}
\define \Hom{\operatorname{Hom}}

\define \cl{\operatorname{cl}}
\define \eps{\varepsilon}
\topmatter
\title
On symplectic manifolds with aspherical symplectic form \endtitle
\author
Yuli Rudyak and Aleksy Tralle
\endauthor
\address
FB6/Mathematik, Universit\"at Siegen, 57068 Siegen, Germany
\endaddress
\vskip6pt
\email
rudyak\@mathematik.uni-siegen.de, 
july\@mathi.uni-heidelberg.de
\endemail
\vskip6pt
\address
 University of Olsztyn, 10561 Olsztyn, Poland 
\endaddress
\vskip6pt
\email
tralle\@tufi.wsp.olsztyn.pl
\endemail
\vskip6pt
\curraddr
Vivatsgasse 7, 53111 Bonn, Germany
\endcurraddr
\email
 tralle\@tufi.wsp.olsztyn.pl, 
tralle\@mpim-bonn.mpg.de
\endemail
\vskip6pt
\date
July 1999
\enddate
\abstract
We consider closed symplectically aspherical manifolds, i.e. closed symplectic manifolds $(M,\omega)$ satisfying the condition $[\omega]|_{\pi_2M}=0$. Rudyak and Oprea [RO] remarked that such manifolds have nice and controllable homotopy properties. Now it is clear that these properties are mostly determined by the fact that the strict category weight of $[\omega]$ equals 2. We apply the theory of strict category weight to the problem of estimating the number of closed orbits of charged particles in symplectic magnetic fields. In case of symplectically aspherical manifolds our theory enables us to improve some known estimations.
\endabstract
\subjclass Primary 55M30, Secondary 58F05 
\endsubjclass
\endtopmatter
\head Introduction
\endhead
A symplectically aspherical manifold is a symplectic manifold $(M,\omega)$ such that  
$$
\int_{S^2}f^*\omega=0
$$
for every map $f: S^2 \to M$. Such manifolds have explicitly appeared in papers of Floer \cite {F} and Hofer \cite{H} in context of Lagrangian intersection theory. These authors used an analytical advantage of symplectic asphericity: namely, it excludes the appearance of non-trivial pseudo-holomorphic spheres in $M$. 

\m However, it turned out that symplectically aspherical manifolds have nice and controllable {\it homotopy} properties. Now it is clear that most of these properties are based on the fact that strict category weight of $[\omega]$ equals 2 (Lemma 2.1), here $[\omega]$ is the de Rham cohomology class of the symplectic form $\omega$. In \cite{RO} a weak version of this result (but with the same proof!) was used in order to prove the equality $\cat M=\dim 
M$ for every closed symplectically aspherical manifold $M$. It was an important ingredient in proving of the Arnold conjecture about symplectic fixed points \cite{R2, RO}.   

\m In this paper we continue to study symplectically aspherical manifolds and, in particular, apply the theory of strict category weight to a special problem in symplectic topology. We show that additional topological arguments enable us to estimate the number of closed orbits of charged particles in {\it symplectic magnetic fields}. This subject has been recently of substantial interest \cite{GK, K}. In a number of recent papers several estimates appeared and it was noted that these are probably not best possible. It follows from the results of this paper that at least in case of {\it symplectically aspherical} closed manifolds the corresponding estimates can indeed be sharpened (Theorem 3.4).
\m
Notice that we don't prove any new analytic results, we only show that our topological observations do help to improve the existing estimates \cite{K} in case of symplectically aspherical manifolds. We think that the technique of category weight potentially may have many applications to problems in symplectic topology, and, in a broader sense, to other non-linear analytic problems.

\m{\bf Acknowledgment}. This article was written when both authors were visiting Max Planck Institut f\"ur Mathematik in Bonn. They are grateful to the Institute and to Professor Yuri Manin for the support of this work. Also, the research of the second author was partially financed by the Polish Research Committee (KBN). 
\vskip6pt
We express our sincere thanks to Ely Kerman and Serge Tabachnikov for useful discussions.
 
\head 1. Topological set-up \endhead
{\mytitle 1.1. Definition {\rm ([LS], [Fox])}.} Given a map $\varphi: A \to X$, we say that a subset $U$ of $A$ is {\it $\varphi$-categorical} if it is open (in $A$) and $\varphi|U$ is null-homotopic. We define the {\it Lusternik--Schnirelmann category} $\cat \varphi $ of $\varphi$ as follows:
$$
\cat \varphi=\min\{k\bigm|A=U_1\cup \cdots \cup U_{k+1} \text{ where each $U_i$ is $\varphi$-categorical}\}.
$$
 Furthermore, we define the {\it Lusternik--Schnirelmann category} $\cat X$ of a space $X$ by setting  $\cat X:= \cat 1_X$.
\proclaim{1.2. Proposition {\rm ([BG])}} {\rm (i)}  For every diagram $A @>\varphi>> Y @>f>> X$ we  have $\cat f\varphi \leq \min\{\cat \varphi, \cat f\}$ whenever $\cat f$ and $\cat \varphi$ are defined. In  particular, $\cat f\leq \min \{\cat X, \cat Y\}$. \par 
{\rm (ii)} If $\varphi \simeq \psi: A \to X$ then $\cat \varphi =\cat \psi$. \par 
{\rm (iii)} If $h: Y \to X$ is a homotopy equivalence then $\cat \varphi=\cat h \varphi$ for every $\varphi: A \to X$.
\qed
\endproclaim
\m Given a path connected pointed space $X$, let 
$$
\eps: S\Omega X \to X
\tag{1.3}
$$ 
be the map adjoint to $1_{\Omega X}$, where $\Omega X$ is the loop space of $X$ and $S$ denotes the (reduced) suspension, see e.g. [Sw].
\proclaim{1.4. Theorem {\rm ([Sv, Theorems 3, 19$'$ and 21])}} Let $\varphi: A \to X$ be a map of connected Hausdorff paracompact spaces. Then $\cat \varphi<2$ iff there is a map $\psi: A \to S\Omega X$ such that $\eps\psi \simeq \varphi$. 
\qed
\endproclaim
{\mytitle 1.5. Definition {\rm ([R1])}.} Let $X$ be a Hausdorff paracompact space, and let $u \in H^q(X;G)$ be an arbitrary element. We define the {\it strict category weight} of $u$ (denoted by $\swgt u$) by setting
$$
\swgt u=\sup\{k\bigm|\varphi^*u=0 \text{ for every map $\varphi: A \to X$ with } \cat \varphi<k\}
$$
where $A$ runs over all Hausdorff paracompact spaces.
\m We use the term ``strict category weight'', since the term ``category weight'' is already used (introduced) by Fadell--Husseini~[FH]. Concerning the relation between category weight and strict category weight, see~[R1].
\m We remark that $\swgt u={\infty}$ if $u=0$. 
\proclaim{1.6. Theorem} Let $X$ be a path connected Hausdorff paracompact space, and let $u\in H^*(X)$ be an arbitrary cohomology class. Then the following holds: \par 
{\rm (i)} if $\deg u>0$ then $\swgt u \geq 1$;\par 
{\rm (ii)} for every map $f: Y \to X$ we have $\cat f \geq \swgt u$ provided that $f^*u \neq 0$; \par
{\rm (iii)} $\cat X \geq \swgt u$ provided that $u\ne 0$;\par
{\rm (iv)} for every map $f: Y \to X$ we have $\swgt f^*u \geq \swgt u$; \par
{\rm (v)} for every $v\in H^*(X)$ we have $\swgt(uv)\geq \swgt u + \swgt v$;\par
{\rm (vi)} $\swgt u\geq 2$ iff $\eps^*u=0$ where $\eps$ is as in $(1.3)$.
\endproclaim
\p (i),(ii) This follows from the definition of $\swgt$.\par
(iii) This follows from (ii).\par
(iv) This follows from 1.2(i).\par
\par (v) Let $\swgt u=k,\ \swgt v=l$ with $k,l<\infty$. Given $f: A \to X$ with $\cat f < k+l$, we prove that $f^*(uv)=0$. Indeed, $\cat f<k+l$, and so $A=A_1\cup \cdots\cup A_{k+l}$ where each $A_i$ is open in $A$ and $f|A_i$ is null-homotopic. Set $B=A_1\cup \cdots\cup A_k$ and $C=A_{k+1}\cup \cdots \cup A_{k+l}$. Then $\cat f|B<k$ and $\cat f|C<l$. Hence $f^*u|B=0=f^*v|C$, and so $f^*(uv)|(B\cup C)=0$. i.e., $f^*(uv)=0$.\par
The case of infinite category weight is left to the reader.\par
(vi) This follows from 1.4.
\qed
\head{2. Strict category weight and symplectically aspherical manifolds}\endhead

Given a space $X$ and an element $u\in H^2(X;G)$, the notation $u|_{\pi_2(X)}=0$ means that
$$
\langle u, h(a)\rangle =0 \text{ for every }a\in \pi_2(X)
$$
where $h: \pi_2(X) \to H_2(X)$ is the Hurewicz homomorphism and $\langle-.-\rangle$ is the Kronecker pairing. 
\m Notice that a symplectic manifold $(M,\omega)$ is symplectically aspherical if and only if $[\omega]|_{\pi_2(M)}=0$.
\proclaim{2.1. Lemma {\rm (cf. [RO])}} Let $X$ be a finite $CW$-space, and let $u\in H^*(X;\Bbb R)$ be a cohomology class such that $u|_{\pi_2(X)}=0$. Then the following holds:\par
{\rm(i)} for every map $f: X \to K(\pi_1(X),1)$ which induces an isomorphism of fundamental groups,
$$
u\in \IM\{f^*: H^2(K(\pi_1(X);\Bbb R) \to H^2(X;\Bbb R)\};
$$
\par
{\rm (ii)} $\swgt u \geq 2$. 
\endproclaim
\p (i) Let $K:=K(\pi_1(X),1)$, and let $f: X \to K$ be a map such that $f_*:\pi_1(X) \to \pi_1(K)$ is an isomorphism. Recall that for every space $Y$ the evaluation homomorphism 
$$
e: H^*(Y;\Bbb R) \to \Hom(H_*(Y),\Bbb R)
$$
is an isomorphism. Furthermore, there is the Hopf exact sequence
$$
\pi_2(X) @>h>> H_2(X) \to H_2(K)  \to 0.
$$
So, we have the following commutative diagram with the exact row:
$$
\CD
 H^2(K;\Bbb R) @>f^*>>  H^2(X;\Bbb R)\\
 @V\cong VV @VeV \cong V @.\\
\Hom(H_2(K),\Bbb R) @>\Hom(f_*,1)>> \Hom(H_2(X),\Bbb R)@>\Hom(h,1)>> \Hom(\pi_2(X),\Bbb R)\\
\endCD
$$
\m Since $\Hom(h,1)(e^*u)=0$, we conclude that $u\in \IM f^*$.
\m (ii) Consider the commutative diagram
$$
\CD
S\Omega X @>S\Omega f >> S\Omega K\\
@V\eps_XVV @VV\eps_K V\\
X @>f>>K
\endCD
$$
with $\eps_X, \eps _K$ as in (1.3). Because of (i), $u=f^*v$ for some $v\in H^2(K;\Bbb R)$. Furthermore, $\Omega K$ is homotopy equivalent to the discrete space $\pi$, and so $S\Omega K$ is homotopy equivalent to a wedge of circles. Hence, $H^2(K;\Bbb R)=0$, and, in particular, $\eps^*_Kv=0$. Hence,
$$
\eps^*_Xu = \eps^*_Xf^*v=(S\Omega f)^*(\eps_Kv)=0
$$
Thus, by 1.6(vi), $\swgt u\geq 2$.
\qed
\proclaim{2.2. Corollary} Let $(M,\omega)$ be a symplectic manifold, and let $K$ denote the space $K(\pi_1(M),1)$. The following three conditions are equivalent:\par
{\rm(i)} $(M,\omega)$ is symplectically aspherical;\par
{\rm(ii)} there exists a map $f: M \to K$ which induces an isomorphism of fundamental groups and such that
$$
[\omega]\in \IM\{f^*: H^2(K;\Bbb R) \to H^2(M;\Bbb R)\};
$$\par
{\rm(iii)}  there exists a map $f: M \to K$ such that
$$
[\omega]\in \IM\{f^*: H^2(K;\Bbb R) \to H^2(M;\Bbb R)\}
$$
\endproclaim
\p (i) $\Rightarrow$ (ii). This follows from 2.1. \par
(ii) $\Rightarrow$ (iii). Trivial.\par
(iii) $\Rightarrow$ (i). If $[\omega]=f^*u$ for some $u\in H^2(K;\Bbb R)$ then $u|_{\pi_2(K)}=0$ since $\pi_2(K)=0$, and thus $[\omega]|_{\pi_2(M)}=0$.
\qed

\proclaim{2.3. Theorem} Let $(M^{2m},\omega)$ be a closed symplectically aspherical manifold, and let
$$
\Bbb CP^n @>i>>E @>p>> M
$$
be a $\Bbb CP^n$-fibration over $M$. Then $\cat E \geq 2m+n$.
\endproclaim
\p Let $u\in H^2(\Bbb CP^n;\Bbb R)$ be a non-zero element.  It is well-known that, for the fibration on hand, its Leray--Serre spectral sequence with real coefficients collapses, [T]. So, there is an element $v\in H^2(E;\Bbb R)$ such that $i^*v=u$. Moreover, according to the Leray--Hirsch Theorem, see e.g. \cite{Sw, 15.47}
$$
(p^*[\omega])^m v^n \neq 0\in H^*(E;\Bbb R).
$$ 
 Now, by 2.1, $\swgt[\omega]\geq 2$, and so, by 1.6(ii,iv,v), 
$$
\cat E\geq \swgt((p^*[\omega])^m v^n)\geq \swgt p^*\left([\omega]^m\right)+\swgt v^n\geq \swgt [\omega]^m +n \geq 2m+n
$$
as claimed.
\qed

\head 3. Application to particle dynamics in magnetic field \endhead

In this section we show how one can apply the previous topological observations to estimation of number of closed trajectories of a charged particle in a symplectic magnetic field. Mathematically, this problem states as follows. 
\m Let $(M,\omega)$ be a closed symplectic manifold, let $p: T^*M \to M$ be the cotangent bundle, and let $\varphi$ be the standard 1-form on $T^*M$. Consider the symplectic form $\Omega:=d\varphi+p^*\omega$ on $T^*M$. Given a Hamiltonian $H$ on $T^*M$, consider the dynamics
$$
i_X\Omega = dH.
\tag{3.1}
$$
The integral curves of the vector field $X$ are the trajectories of a charged particle in a magnetic field with the Hamiltonian $H$.

\m Choose a Riemannian metric $g$ and an almost complex structure $J$ on $M$ such that 
$$
\omega(Y,Z)=g(Y,JZ)
$$
for every vector fields $Y,Z$ on $M$ and consider the so-called metric Hamiltonian 
$$
H: T^*M \to R,\quad H(\xi)=\langle \xi,\xi\rangle_g
$$
for every covector $\xi$ on $M$.  
Kerman~[K] proved that, for every $\eps>0$ small enough, the energy level $H^{-1}(\eps)$ contains at least $m+\cl(M)$ closed trajectories of the vector field $X$  as in (3.1); here $\cl$ denotes the cup-length. In greater detail, Kerman proved the following
\proclaim{3.2. Theorem {\rm([K])}} For any metric Hamiltonian $H$ and every $\eps$ small enough, there is a fiber bundle $p: \Sigma \to M$ with the following properties:
\roster
\item $\Sigma\subset H^{-1}(\eps)$;
\item $p: \Sigma \to M$ is a locally trivial bundle with the fiber $S^{2m-1}$;
\item there is a free fiberwise $S^1$-action on $\Sigma$, and the restriction of this action to any fiber coincides with the standard $S^1$-action on $S^{2m-1}$;
\item the number of periodic trajectories of the vector field $X$ is at least $1+\cat (\Sigma/S^1)$. \footnotemark \footnotetext {Notice that the Lusternik--Schnirelmann category considered in our paper is one less then the Lusternik--Schnirelmann category used in [K]}
\endroster 
\endproclaim
Notice that the projection $p: \Sigma \to M$ yields a locally trivial bundle
$$
q: \Sigma/S^1 \to M
$$
with the fiber $S^{2m-1}/S^1=\Bbb CP^{m-1}$. Because of this, Kerman obtains the estimation
$$
\cat (\Sigma/S^1)\geq \cl(\Sigma/S^1)\geq \cl M +\cl(\Bbb CP^{m-1})=\cl M +m-1
$$
and so 
$$
\text {\it the number of closed trajectories of $X$ is at least $m+\cl M $.}
\tag{3.3}
$$
\m We remark that this estimation can be improved if the symplectic form $\omega$ on $M$ is aspherical. 

\proclaim{3.4. Theorem} Let $(M^{2m},\omega)$ be a closed symplectic manifold with $[\omega]|_{\pi_2(M)}=0$, and let $H$ be a metric Hamiltonian on $T^*M$. Then, for every $\eps$ small enough, the set $H^{-1}(\eps)$ contains at least $3m$ closed trajectories of the vector field $X$ as in $(3.1)$.
\endproclaim
\p Since $\Sigma/S^1$ is fibered over $M$ with the fiber $\Bbb CP^{ m-1}$, we conclude that, by Theorem 2.3, $\cat(\Sigma/S^1)\geq 3m-1$. Now the result follows from 3.2(4). 
\qed
{\mytitle 3.5. Remarks.} 1. Using the language of mathematical physics,  3.4 can be reformulated as follows. {\it A charged particle on a symplectically aspherical manifold, under the influence of the symplectic magnetic field, has at least $3m$ closed trajectories on any sufficiently low kinetic energy level.}
\par 2. Clearly, $\cl(M)\leq 2m$ for every connected $2m$-dimensional manifold $M$. Furthermore, there are many examples of closed symplectic manifolds with aspherical symplectic form and such that $\cl M <2m$, see Proposition 4.3 below. So, for closed aspherical symplectic manifolds the estimation from 3.4 is stronger than this one from (3.3).   

\head 4. Examples
\endhead

In this section we present examples of symplectically aspherical closed manifolds. These examples give the reader a picture of the class of manifolds covered by the results of the present paper.

\m To start with, notice that there are two constructions which preserve symplectic asphericity. First, the Cartesian product of two symplectically aspherical manifolds is symplectically aspherical. Second, if $f: (M,\omega) \to (N,\sigma)$ is a map of symplectic manifolds with $f^*\sigma=\omega$, then $(M,\omega)$ is symplectically aspherical provided $(N,\sigma)$ is. This can be proved directly or deduced from 2.2.
\par As a special case, consider a symplectically aspherical maniflod $(M,\omega)$ and a branched covering $p: \T M \to M$ such that the branch locus is a symplectic submanifold of codimension 2. Then $\T M$ possesses a symplectic form $\T \omega$ with $[\T \omega]=p^*[\omega]$, [G]. Thus, $(\T M, \T \omega)$ is symplectically aspherical. 
 
\m Notice that if a manifold, say, $M$ is covered by $\Bbb R^n$, then it is  necessarily aspherical, i.e. $\pi_iM=0$ for $i>1$. Therefore, for any symplectic  manifold covered by $\Bbb R^n$ we have $\pi_2M=0$ and, in particular, the condition $[\omega]|_{\pi_2M}=0$ is trivially satisfied. So, taking quotients of simply connected solvable Lie groups by lattices (if there are any), we get many examples of aspherical closed manifolds.

\example{4.1. Example {\rm (nilmanifolds and solvmanifolds)}} Let $N$ be any simply connected nilpotent Lie group with the Lie algebra $\frak n$ satisfying the following condition: there exists a basis $e_1,..., e_n$ of $\frak n$ such that all the structural constants $c_{ij}^k$ of $\frak n$ with respect to this basis are {\it rational numbers}. By the Malcev Theorem, $N$ necessarily admits a lattice, i.e. a discrete co-compact subgroup $\Gamma$. Hence we can form a closed aspherical nilmanifold $N/\Gamma$. One can easily get a characterization of {\it symplecticness} of such manifolds. It is known that each cohomology class in $H^*(N/\Gamma)$ contains a {\it homogeneous} representative, i.e. a differential form whose pullback to $N$ is an invariant differential form on Lie group $N$, see \cite{TO, Ch. 2, Theorem 1.3}. Hence, every cohomologically symplectic manifold of the form $N/\Gamma$ is symplectic.  (Recall that a closed connected smooth manifold $M^{2n}$ is called {\it !
cohomologically symplectic} if it possesses a 2-form $\omega$ such that $[\omega]^n\ne 0$.) One can check the cohomological symplecticness by calculating the Chevalley-Eilenberg complex for $\frak n$. So, we have the following fact.
\endexample

\proclaim{4.2. Proposition} Any closed cohomologically symplectic nilmanifold $M=N/\Gamma$ is simplectically aspherical.
\qed
\endproclaim

Now we present an explicit example of a symplectically aspherical closed 4-dimensional manifold $K$ for which $\cl K =3$ and $\cat K=4$. Consider the {\it Kodaira-Thurston manifold} $K=(N_3/\Gamma)\times S^1$ defined as follows. Let $N_3$ stand for the 3-dimensional Heisenberg group of all matrices of the form
$$
\pmatrix
1 & a & b \\
0 & 1 & c\\
0 & 0 & 1
\endpmatrix
$$
where $a, b, c\in\Bbb R$ and $\Gamma$ denotes the subgroup of $N_3$ consisting of all matrices with integer entries. The cohomology of $K$ is easy to calculate as was indicated below (see also \cite{TO}): the Chevalley--Eilenberg complex of $K$ is given by the formulae
$$(\Lambda(\alpha,\beta,\gamma,\delta),d)$$
$$d\alpha=d\beta=d\delta=0,\quad d\gamma=\alpha\wedge\beta$$
where $\Lambda$ denotes the exterior algebra generated by $\alpha,\beta,\gamma,\delta$ and these generators have degree 1. Hence the vector space $H^1(K;\Bbb R)\simeq \Bbb R^3$ is generated by $[\alpha],[\beta],[\delta]$ and, obviously, the length of the longest product in $H^*(K;\Bbb R)$ is $\leq 3$. In fact, $\cl(K)=3$, the element $[\alpha][\delta][\beta\gamma]$ gives us the product of length 3. Notice that $K$ is symplectic, since the cohomology class $[\alpha\wedge \delta+\beta\wedge\gamma]$ is symplectic. However, $\text{cat}\,(K)=4$, since $\text{cat}\,V=\dim\,V$ for any closed aspherical manifold $V$ (see \cite{EG}). Thus,
$$
3=\cl K<\cat K=4.
$$
 
\m The previous example can be obviously generalized to any symplectic nilmanifold. Indeed, the Chevalley--Eilenberg complex of any nilmanifold $N/\Gamma$ has the form
$$(\Lambda(x_1,...,x_n),d)$$
with $n=\dim\,(N/\Gamma)=\dim\,N$ and the degrees of generators $x_i$ being 1. If $N/\Gamma$ is not a torus (i.e.,    if $\Gamma$ is not Abelian), then $ \Ker d \ne 0$,  and so $dx_i\not=0$ for at least one $x_i$. Hence, any cup-product of any $n$ cocycles is a boundary, and we obtain the following propsition.

\proclaim{4.3 Proposition} For any cohomologically symplectic non-toral nilmanifold $N/\Gamma$ 
$$
\cl(N/\Gamma)<n=\dim\,(N/\Gamma)=\cat(N/\Gamma).\qed
$$
\endproclaim

\m In the same fashion one can analyze {\it solvmanifolds}. Here, however, it is not so easy to prove the existence of lattices and symplectic structures (cf. \cite{TO}). Nevertheless, there is a certain class of solvmanifolds where the answer is essentially the same as that for nilmanifolds.  We say that a simply connected solvable Lie group $G$ is {\it completely solvable}, if for the Lie algebra $\frak g$ of $G$ all the operators 
$$
\text{ad}\,V: \frak g\to \frak g,\quad V\in \frak g
$$ 
have only real eigenvalues. Now, the following proposition can be proved similarly to 4.2.

\proclaim{4.4. Proposition}  Let $\Gamma$ be a co-compact discrete subgroup of a completely solvable Lie group $G$. If the manifold $M=G/\Gamma$ is cohomologically symplectic, then it is symplectically aspherical.
\qed
\endproclaim

Notice that there are also examples of symplectically aspherical closed manifolds with non-trivial $\pi_2$, but the corresponding explicit constructions are rather delicate and have very recently appeared \cite{G}. The general description of the whole class of symplectically aspherical manifolds is given in Corollary 2.2. However, it seems rather difficult to use it to construct explicit examples. Recently Gompf produced such examples as branched coverings.

\proclaim{4.5. Theorem {\rm ([G])}} Let $p:\T X \to X$ be a $d$-fold branched covering of an orientable 4-manifold $X$  such that the  branch locus $B\subset X$ is obtained from a generically immersed surface $B^*\subset X$ by smoothing all double points. If $B^*$ has $k\geq 1$ double points then $\pi_2(\T X)\otimes \Bbb R\neq 0$.
\qed
\endproclaim
\m So, if one takes any 4-dimensional closed symplectically aspherical manifold $(X,\omega)$ (e.g., the Kodaira--Thurston manifold) and a branched covering  $p: \T X \to X$ as in 4.5, then $(\T X, \T \omega)$ (described in the beginning of the section) is a simplectically aspherical manifold with $\pi_2(\T X)\ne 0$.

\example{4.6. Example} We can get more aspherical closed manifolds by considering lattices in {\it semisimple} Lie groups of non-compact type. It is known \cite{FJ} that any simply connected semisimple Lie group contains a torsion-free, discrete, co-compact subgroup. Consequently, any symmetric space of non-compact type is the universal covering space of a locally symmetric space $M$. Hence, $M$ is a closed and aspherical manifold.
\vskip6pt
In the same fashion, we can consider {\it hyperbolic manifolds} (examples: Riemannian surfaces).
\endexample

\enddocument
\end